\documentclass[12pt]{amsart}
\usepackage{latexsym, amssymb, amsmath, amscd}
 \usepackage{eucal}

\theoremstyle{plain}
    \newtheorem{rema}{Remark}[section]

   \newtheorem{theo}[rema]{Theorem}
 \newtheorem{conj}[rema]{Conjecture}
   
    \newtheorem{lemma}[rema]{Lemma}
    \newtheorem{corol}[rema]{Corollary}
     \newtheorem{exam}[rema]{Example}
  \newtheorem{rmk}[rema]{Remark}

	\newcommand{\nno}{\nonumber}

	\newcommand{\p}{\partial}

 \newcommand{\pf}{{\it Proof:}\hspace{2ex}}
 
 \newcommand{\epfv}{\hspace{1em}$\Box$\vspace{1em}}

\newcommand{\bZ}{{\mathbb Z}}
\newcommand{\bR}{{\mathbb R}}
\newcommand{\bQ}{{\mathbb Q}}
\newcommand{\bN}{{\mathbb N}}

\newcommand{\poly}{ \mbox{\rm Poly\,}}
\newcommand{\supp}{ \text{\rm Supp\,} }

\newcommand{\bC}{\mathbb C}

\newcommand{\cz}{ \bC[z] }
\newcommand{\czz}{ \bC[z^{-1}, z] }
\newcommand{\qrt}{ \bR_{\ge 0}^n}

\newcommand{\lin}{1\le i\le n}

\renewcommand{\theequation}{\thesection.\arabic{equation}}
\renewcommand{\therema}{\thesection.\arabic{rema}}
\newcommand{\EAn}{\end{align*}}

\title[Some Results on the Vanishing Conjecture]
{Some Results on the Vanishing Conjecture of Differential Operators with Constant Coefficients}

    \author{Arno van den Essen, Roel Willems and Wenhua Zhao}      
   \date{\today}
    
\begin{document}

\begin{abstract}
In this paper we prove four cases of the vanishing conjecture of differential operators with constant coefficients and also a conjecture on the Laurent polynomials with no holomorphic parts, which were proposed in \cite{GVC} by the third named  author. We also give two examples to show that the generalizations of both the vanishing conjecture and 
the Duistermaat-van der Kallen theorem \cite{DK} 
to Laurent formal power series do not hold in general.   
\end{abstract}

\keywords{The vanishing conjecture of differential operators 
with constant coefficients, polytopes of Laurent polynomials, 
the Jacobian conjecture}
   
\subjclass[2000]{33C45, 32W99, 14R15}

\thanks{The second named author is supported by the Netherlands Organization of Scientific Research (NWO). The third named  author has been partially supported by NSA Grant R1-07-0053}

 \bibliographystyle{alpha}
    \maketitle


\renewcommand{\theequation}{\thesection.\arabic{equation}}
\renewcommand{\therema}{\thesection.\arabic{rema}}
\setcounter{equation}{0}
\setcounter{rema}{0}
\setcounter{section}{0}

\section{\bf Introduction}\label{S1}

Let $z=(z_1, z_2, \dots, z_n)$ be $n$ commutative 
free variables and $\cz$ (resp.\,$\czz$) the algebra of polynomials (resp.\,Laurent polynomials) in $z$ over 
$\bC$. For any $\lin$, set $\p_i\!:=\p/\p z_i$,  
and $\p\!:=(\p_1, \p_2, ..., \p_n)$.

We say a differential operator $\Lambda$ of $\cz$ 
is a 
{\it differential operator with constant coefficients}  
if $\Lambda=h(\p)$ for some polynomial 
$h(\xi)\in \bC[\xi]$, where 
$\xi=(\xi_1, \xi_2, \dots, \xi_n)$ are 
another $n$ commutative variables which also 
commute with $z$. For convenience, 
we will denote the polynomial 
$h(\p)$ by $\Lambda(\p)$ 
and simply write $\Lambda=\Lambda(\p)$.

In this paper, we will prove four cases of the following 
{\it vanishing conjecture} of differential operators 
with constant coefficients, which was proposed 
by the third named author in 
\cite{GVC}.  

\begin{conj}\label{VC}
Let $P(z)\in \cz$ and $\Lambda=\Lambda(\p)$ 
for some $\Lambda(\xi)\in\bC[\xi]$. Assume 
that $\Lambda^m (P^m)=0$ for any $m\ge 1$. Then 
$\Lambda^m (P^{m+1})=0$ 
when $m\gg0$.
\end{conj}

Actually, all the cases of the conjecture above that 
we will prove in this paper also hold in 
the following more general form.    

\begin{conj}\label{GVC}
Let $P(z)\in \cz$ and $\Lambda=\Lambda(\p)$ 
for some $\Lambda(\xi)\in\bC[\xi]$. Assume 
that $\Lambda^m (P^m) =0$ for any $m\ge 1$. Then, 
for any $g(z)\in \cz$, we have $\Lambda^m(P^mg)=0$ 
when $m\gg0$.
\end{conj}

Note that, Conjecture \ref{VC} is just the special 
case of Conjecture \ref{GVC} with $g(z)=P(z)$.
Also, when $\Lambda$ is a homogeneous 
differential operator of order $2$ with constant 
coefficients, these two conjectures are 
actually equivalent (See \cite{EZ} 
and \cite{GVC}).

Note also that Conjecture \ref{GVC} has recently been   generalized by the third named author \cite{IC} 
to the so-called  
{\it image conjecture} of commuting differential operators of 
order one with constant leading coefficients. Actually, 
Conjecture \ref{GVC} is equivalent to the 
image conjecture of the commuting 
differential operators $\xi_i-\p_i$ 
$(\lin)$ of the polynomial algebra 
$\bC[\xi, z]$ for the separable polynomial  
$\Lambda(\xi)P(z)\in \bC[\xi, z]$. 
For more details, see \cite{IC}.
 
The main motivation behind Conjecture 
\ref{VC} is its connection with the well-known 
{\it Jacobian conjecture} proposed by O. H. Keller 
\cite{K} in $1939$ (See also \cite{BCW} and \cite{E}).
The connection is given by the following theorem 
proved in \cite{HNP}.

\begin{theo}\label{JC-VC}
Let $\Delta_n \!:=\sum_{i=1}^n \p_i^2$ be 
the Laplace operator of the polynomial 
algebra $\cz$. Then, the Jacobian conjecture holds for all 
$n\ge 1$ iff Conjecture \ref{VC} with 
$\Lambda=\Delta_n$ holds for all $n\ge 1$.
\end{theo} 

It has also been shown in \cite{GVC} 
that one may replace the Laplace operators 
$\Delta_n$ in the theorem above by any 
sequence 
$\Lambda_n=\Lambda_n(\p)$ $(n\ge 1)$ 
of differential operators 
with $\Lambda_n(\xi)$ homogeneous of degree 
$2$ whose ranks go to $\infty$ as 
$n\to \infty$.

The proof of Theorem \ref{JC-VC} is based on 
the remarkable symmetric reduction achieved 
independently by M. de Bondt and the first named author 
\cite{BE1} and G. Meng \cite{Me}. It also depends 
on some results obtained in \cite{BurgersEq} 
on a deformation of polynomial maps.  

Currently, there are only a few cases of 
Conjecture \ref{VC} that are known. The best results so 
far come from M. de Bondt and the first named author's 
results \cite{BE2} \cite{BE3} on symmetric 
polynomial maps via the equivalence 
obtained by the third named author 
in \cite{HNP} and \cite{GVC}. 
The results state that 
Conjecture \ref{VC} holds for 
homogeneous quadratic differential operators 
$\Lambda=\Lambda(\p)$ if either $n\le 4$ 
(for any $P(z)\in \cz$) or 
$n\le 5$ with $P(z)$ homogeneous. 
The case when $\Lambda(\xi)$ is an (integral) 
power of a homogeneous linear polynomial of $\xi$ 
is an easy exercise (See also \cite{GVC}). 

In this paper, we will prove four more cases 
of Conjecture \ref{GVC} and also a conjecture 
proposed in \cite{GVC} by the third named author 
on the Laurent polynomials with no holomorphic parts 
(See Theorem \ref{H-Case}).

First, in Section \ref{S2}, we use a fundamental 
theorem, Theorem \ref{ODE}, in ODE to 
show that Conjecture \ref{GVC} holds for the 
one variable case (See Theorem \ref{OneVariableVC}). 
Actually, in this case Conjecture \ref{GVC} 
even holds for all formal power series 
$P(z)$ and polynomials $g(z)\in \cz$.

In Section \ref{S3}, we assume $n=2$ and 
show in Theorem \ref{1+P-Case} 
that Conjecture \ref{GVC} holds for any 
differential operator $\Lambda$ of 
the form $\p_x-\Phi(\p_y)$, 
where, for convenience, 
in this section we use $(x, y)$ instead 
of $(z_1, z_2)$ to denote 
two free commutative variables and 
$\Phi(\cdot)$ to denote any polynomial 
in one variable. We also give 
an example, Example \ref{Counter-GVC}, 
to show that the generalizations 
of both Conjectures \ref{VC} and \ref{GVC} 
to formal power series do not hold 
in general.   

In Section \ref{S4}, we first recall and prove  
some results on rational polytopes, 
polytopes with all its vertices having 
rational coordinates. We then use 
the remarkable 
Duistermaat-van der Kallen theorem \cite{DK} 
(See Theorem \ref{ThmDK}) 
to show what we call the {\it density theorem}  
of polytopes of Laurent polynomials  
(See Theorem \ref{DensityThm}). We will also 
show in Lemma \ref{Pnot>L} that Conjecture 
\ref{GVC} holds when the polytope 
$\poly(P)-\poly(\Lambda)$ has 
no intersection points with 
$(\bR^{\ge 0})^{\times n}$. At the end of this section, 
we give an example, Example \ref{Counter-ThmDK}, 
to show that the Duistermaat-van der Kallen theorem 
can not be generalized to the setting of 
Laurent formal power series.    

In Section \ref{S5}, we first give a proof for  
a conjecture proposed in \cite{GVC} 
on Laurent polynomials with 
no holomorphic parts 
(See Theorem \ref{H-Case}). We then show 
in Corollary \ref{M-Case} that Conjecture 
\ref{GVC} holds when either $\Lambda(\xi)$ or $P(z)$ 
is a monomial of $\xi$ or $z$, respectively. 

In Section \ref{S6}, we show in Theorem \ref{MM-Case} and Corollary \ref{MM-Case2} that Conjecture 
\ref{GVC} holds when either $\Lambda(\xi)$ or $P(z)$ 
is a linear combination of two monomials of $\xi$ 
or $z$, respectively, with different degrees. \\

{\bf Acknowledgment}
The authors would like to thank Harm Derksen, 
Wilberd van der Kallen and Han Peters for personal communications. The authors would also like to thank the anonymous referee for some valuable suggestions.

\renewcommand{\theequation}{\thesection.\arabic{equation}}
\renewcommand{\therema}{\thesection.\arabic{rema}}
\setcounter{equation}{0}
\setcounter{rema}{0}

\section{\bf Proof of the Vanishing Conjecture for One Variable Case}
\label{S2}

In this section, we consider Conjecture \ref{GVC} 
for the one-variable case and show that it does 
hold even for formal power series $P(z)$ 
(See Theorem \ref{OneVariableVC}).

Throughout this section we assume that $n=1$ and let $z$ denote 
a single free variable. For convenience, 
we also set $D\!:=d/dz$. 
A different operator $\Lambda=\Lambda(D)$ is also fixed, 
where $\Lambda(\xi)$ denotes any non-zero polynomial 
in one variable.

With the notation above, let us first recall 
the following two well-known fundamental 
results from ODE (See \cite{L} or any standard 
text book on ODE). 

\begin{theo}\label{ODE}
Let $\lambda_i$ $(1\le i\le k)$ be the set of all 
distinct roots of the polynomial $\Lambda(\xi)\ne 0$ 
with multiplicity $m_i\ge 1$. 
Then, a formal power series 
$P(z)\in \bC[[z]]$ satisfies the 
differential equation $\Lambda(D)P(z)=0$ iff $P(z)$ 
can be written as a linear combination over $\bC$ of 
$z^je^{\lambda_i z}$ $(1\le i \le k;\, 0\le  j \le m_i-1)$.
\end{theo}

\begin{lemma}\label{Uniqueness}
For any distinct $\lambda_i\in \bC$ 
$(1\le i\le k)$, the formal power series 
$\{ e^{\lambda_i z}\,|\, 1\le i\le k\}$ are 
linearly independent over the rational function 
field $\bC(z)$.
\end{lemma}

In case that a formal power series $P(z)$ can be written (uniquely) as 
\begin{align}\label{Exp-expansion}
P(z)=\sum_{i=1}^k c_i(z)e^{\lambda_i z} 
\end{align}
for some distinct $\lambda_i\in \bC$ 
$(1\le i\le k)$ and non-zero $c_i(z)\in \bC[z]$, 
we call this expression the {\it exponential expansion} 
of $P(z)$.

The main result of this section is the following theorem.

\begin{theo}\label{OneVariableVC}
For any formal power series $P(z)$ with 
$\Lambda^m (P^m)=0$ for any $m\ge 1$, we have 

$(a)$ $P(z)$ must be a polynomial in $z$.

$(b)$ Conjecture \ref{GVC} holds for $\Lambda$, $P(z)$ and 
any $g(z)\in \bC[z]$.
\end{theo}

\pf In the proof below, we will freely use the notation 
fixed above for the differential operator 
$\Lambda=\Lambda(D)$ and the polynomial $\Lambda(\xi)$.

Note first that, if $P(z)=0$, 
there is nothing to prove. So we assume 
$P(z)\ne 0$. 

Since $\Lambda P=0$, by Theorem \ref{ODE} and without 
losing any generality, 
we may write $P(z)$ uniquely as in 
Eq.\,(\ref{Exp-expansion}) with $\lambda_i$ 
$(1\le i\le k)$ being distinct roots of 
$\Lambda(\xi)$ and $c_i(z)\in \cz$ with   
$\deg c_i(z)\le m_i-1$ for any $1\le i\le k$.

Assume that $P(z)$ is not a polynomial. Then, 
by Lemma \ref{Uniqueness}, there exists 
$\lambda_i \ne 0$ for some $1\le i\le k$. 
Identify $\bC$ with $\bR^2$ and let $\Sigma$ 
be the polytope or the convex subset in $\bR^2$ 
generated by $\lambda_i$ $(1\le i\le k)$. Then, 
there exists a vertex of $\Sigma$ 
which is not the origin $0\in \bR^2$. 
Without losing any generality, we assume 
that $\lambda_1$ is such a vertex.

For any $m\ge 1$, from Eq.\,(\ref{Exp-expansion}),
it is easy to see that $P^m(z)$ also 
has an exponential expansion 
in which $e^{m\lambda_1 z}$ appears with 
the nonzero coefficient $c_{\lambda_1}^m(z)$.

On the other hand, since $\Lambda^m (P^m)=0$, 
by Theorem \ref{ODE} and Lemma \ref{Uniqueness},
we know that $m\lambda_1$ must be 
a root of the polynomial 
$\Lambda^m(\xi)$, hence also 
a root of $\Lambda(\xi)$. 
Since $\lambda_1\ne 0$ and 
the statement above holds for any $m\ge 1$, 
we see that $\Lambda(\xi)$ has infinitely 
many distinct roots $m\lambda_1$ $(m\ge 1)$, 
which is impossible. Therefore, $P(z)$ must be 
a polynomial and $(a)$ holds.

To show $(b)$, by $(a)$ and Lemma 
\ref{Uniqueness}, it is easy see that 
none of non-zero roots of $\Lambda(\xi)$ 
can be involved in Eq.\,(\ref{Exp-expansion}).
Hence we have $k=1$, $\lambda_1=0$ and 
$P(z)=c_1(z)$ with the degree 
\begin{align}\label{OneVariableVC-pe2}
d\!:=\deg P(z)=\deg c_1(z) \le m_1-1.
\end{align}
Since $\Lambda(\xi)$ has the root 
$\lambda_1=0$ with multiplicity $m_1\ge 1$,
we may write $\Lambda(\xi)$ as 
$\Lambda(\xi)=\Phi(\xi) \xi^{m_1}$
for some $\Phi(\xi) \in \bC[\xi]$.
Consequently, we have 
$\Lambda=\Lambda(D)=\Phi(D) D^{m_1}$.

Now we fix any $g(z)\in \cz$ with $d'\!:=\deg g(z)\ge 0$.
Then, for any $m\ge 1$, the polynomial 
$D^{m m_1}(P^m g(z))$, if not zero, has the degree
\begin{align}\label{OneVariableVC-pe3}
\deg D^{m m_1}(P^m g(z))=(md+d')-m m_1=d'-(m_1-d)m. 
\end{align}

Note that, by Eq.\,(\ref{OneVariableVC-pe2}), 
we know that $m_1-d\ge 1$. So, for any 
$m > d'/(m_1-d)$, we have $d'-(m_1-d)m<0$.
Furthermore, by Eq.\,(\ref{OneVariableVC-pe3}), 
we have, for any $m > d'/(m_1-d)$,
$D^{m m_1}(P^m g(z))=0$, hence also
\begin{align*}
\Lambda^m (P^m g(z))=\Phi^m D^{m m_1}(P^m g(z))=0.
\end{align*}
Therefore, Conjecture \ref{GVC} holds for 
any $g(z)\in \cz$.
\epfv

Two remarks about Theorem \ref{OneVariableVC} are as follows. 

First, the theorem does not always hold for formal power series 
$g(z)\in \bC[[z]]$. For example, 
let $c \in \bC$ such that $c$ is not a root of 
$\Lambda(\xi)$ and $g(z)=e^{cz}$. Then, By Theorem  
\ref{OneVariableVC}, $(a)$ and Theorem \ref{ODE}, 
it is easy to see that $\Lambda^m (P^m(z)g(z))\ne 0$ 
for any $m\ge 1$.

Second, even though the conjecture \ref{GVC} fails 
for some formal power series $g(z)\in \bC[[z]]$,  
by Theorem \ref{OneVariableVC}, $(a)$, 
it still holds when $g(z)=P(z)$ 
without assuming in advance that 
$g(z)=P(z)$ is a polynomial. 
In other words, Conjecture \ref{VC} 
actually still holds for all formal 
power series $P(z)$.

\renewcommand{\theequation}{\thesection.\arabic{equation}}
\renewcommand{\therema}{\thesection.\arabic{rema}}
\setcounter{equation}{0}
\setcounter{rema}{0}

\section{\bf Proof of the Vanishing Conjecture for the Differential Operator $\Lambda=\p_x-\Phi(\p_y)$}\label{S3} 

Throughout this section, we denote by $(x, y)$ instead of 
$(z_1, z_2)$ two commutative free variables and by $\xi$ another free variable which commutes with $x$ and $y$. Once and for all, we also fix an arbitrary non-zero polynomial 
$\Phi(\xi)$ and write it as 
\begin{align}\label{Def-Phi} 
\Phi(\xi)=q_0+q_1\xi+\cdots +q_k \xi^k
\end{align} 
for some $k\ge 0$ and $q_i\in \bC$ $(0\le i\le k)$. 

We will denote by $o(\Phi(\xi))$ or $o(\Phi)$ the {\it order} of the polynomial $\Phi(\xi)$, i.e.\@ the least integer 
$m\ge 0$ such that $q_m\ne 0$.

In this section, we first give a proof of Conjecture \ref{GVC} for the differential operator $\Lambda=\p_x-\Phi(\p_y)$ (See the theorem below). We then give an example 
(See Example \ref{Counter-GVC}) to 
show that the generalizations of both 
Conjectures \ref{VC} and \ref{GVC} to 
formal power series $P(x, y)\in \bC[[x, y]]$ actually 
do not hold.

The main result of this section 
is the following theorem.

\begin{theo}\label{1+P-Case}
The Conjecture \ref{GVC} holds for the differential operator 
$\Lambda=\p_x-\Phi(\p_y)$ and all polynomials 
$P(x, y)\in \bC[x, y]$.
\end{theo}

Note first that, if $P(x, y)=0$, there is nothing to prove for 
the theorem. So, for the rest of this section, we fix an arbitrary polynomial $P(x, y)\in \bC[x, y]$ and  
assume $P(x, y)\ne 0$.

In order to prove Theorem \ref{1+P-Case}, we first need the following three lemmas. 

\begin{lemma}\label{lemma1}
Let $\Lambda$ and $0\ne P(x,y)\in \bC [x,y]$ as fixed above. Assume that $\Lambda P=0$. Then, $q_0=0$, or equivalently, the order $o(\Phi(\xi))\geq 1$ if $\Phi(\xi)\ne 0$.
\end{lemma}
\pf Assume $q_0\ne 0$. By Eq.\,(\ref{Def-Phi}), we have 
$\Lambda = \p_x - \Phi(\p_y) 
= \p_x -(q_0 + q_1\p_y +\cdots + q_k\p_y^k)$. 
Let $d=\deg P(x,y)$ and $P_d(x,y)$ 
the homogeneous part of $P(x,y)$ of degree $d$.
Note that the highest degree part of $\Lambda P$ is 
$q_0P_d(x,y)$ which is equal to zero since $\Lambda P=0$. 
Since we have assumed $q_0\ne 0$ and $P(x, y)\ne 0$, 
we have $P_d(x,y)=0$, which is a contradiction.
\epfv

\begin{lemma}\label{lemma2}
Let $\Lambda$ and $0\ne P(x,y)\in \bC [x,y]$ as above with 
$\Lambda P=0$. Then, we have $P(x,y)=e^{x\Phi(\p_y)}(f(y))$ 
for some $f(y)\in \bC[y]$.
\end{lemma}
\pf Note first that, 
\begin{align*}
\p_x(e^{-x\Phi(\p_y)}P) = e^{-x\Phi(\p_y)}(\p_x-\Phi(\p_y))P
=e^{-x\Phi(\p_y)}(\Lambda P)=0.
\end{align*}
So $e^{-x\Phi(\p_y)}P$ is independent on $x$. Hence 
$e^{-x\Phi(\p_y)}P=f(y)$ for some $f(y)\in \bC[y]$. 
Now, applying $e^{x\Phi(\p_y)}$ to both sides of the latter equation, we get $P(x,y)=e^{x\Phi(\p_y)}(e^{-x\Phi(\p_y)}P)
=e^{x\Phi(\p_y)}(f(y))$.
\epfv

Note that $e^{x\Phi(\p_y)}(f(y))$ is still a polynomial, because from Lemma \ref{lemma1} it follows that $\Phi=0$ or $o(\Phi)\geq 1$. So, in the first case 
$P(x, y)=f(y)$ and in the latter there exists an $m$ such that 
$\Phi^m(\p_y)(f(y))=0$.

\begin{lemma}\label{lemma3}
Let $\Lambda$ and $P(x,y)$ be as in Lemma \ref{lemma2}. Further  
assume $o(\Phi)\ge 2$ and $\Lambda P=\Lambda^2(P^2)=0$. 
Then, we have $o(\Phi)>\deg f$ and $P(x,y)=f(y)$.
\end{lemma}

\pf First, we view 
$\Lambda ^2(P^2) = 
(\p_x - \Phi(\p_y))^2( (e^{x\Phi(\p_y)}f(y))^2 )$ as a polynomial in $\bC[y][x]$,
and look at its constant term, which is 
\begin{align*}
\Lambda ^2(P^2)_{|x=0} & =  2f(y)\Phi^2(\p_y)f(y) + 
2(\Phi(\p_y)f(y))^2 \\
 & \quad\quad -4\Phi(\p_y)\left(f(y)\Phi(\p_y)f(y)\right) + 
\Phi^2(\p_y)(f^2(y))\\
 & = 0 
\end{align*}

Let $d=\deg f(y)$ and $r=o(\Phi(\xi))\geq 2$, 
and assume that $d\geq r$. 
Write $f(y) = c_0 + c_1y + \cdots + c_dy^d$ and $\Phi(\p_y) = q_r\p_y^r + q_{r+1}\p_y^{r+1} + \cdots + q_{k}\p_y^{k}$ for some 
$k\geq r$.

Now, by looking at the leading coefficient of 
$\Lambda ^2(P^2)_{|x=0}$, we get 
\begin{align*}
0&=2c_d^2q_r^2 v\frac{d!}{(d-r)!} + 2c_d^2q_r^2\left(\frac{d!}{(d-r)!}\right)^2 \\
 &\quad \quad -4c_d^2q_r^2\frac{d!}{(d-r)!}\frac{(2d-r)!}{(2d-2r)!} + c_d^2q_r^2\frac{(2d)!}{(2d-2r)!},
\end{align*}
where $v=0$ if $d<2r$ and $v=\frac{(d-r)!}{(d-2r)!}$ 
if $d\geq 2r$.

Then, by the assumptions $c_d \not= 0$ and $q_r \not= 0$, we further have 
\begin{align}\label{eqn}
 0=\frac{2vd!}{(d-r)!} + 
2\left(\frac{d!}{(d-r)!}\right)^2 
- \, \frac{4d!(2d-r)!}{(d-r)!(2d-2r)!} + 
\frac{(2d)!}{(2d-2r)!}.
\end{align}

But, on the other hand, for any $r\ge 2$, we also have 
\begin{align}
\frac{2vd!}{(d-r)!} + 2\left(\frac{d!}{(d-r)!}\right)^2  
& >  0 \label{ineq-1} \\
- \, \, \frac{4d!(2d-r)!}{(d-r)!(2d-2r)!}+\frac{(2d)!}{(2d-2r)!}  & \geq  0. \label{ineq-2}
\end{align}
The first inequality is obvious. The second inequality holds 
is because it is equivalent to
$\binom{2d}{r} \geq  4\binom{d}{r}$ which  
follows from the general identity 
$\binom{2d}{r}\geq 2^r\binom{d}{r}$ for any $r\geq 0$ 
and the assumption that $r\geq 2$.
Note that the last inequality can be easily verified by the facts that $2d-i\ge 2(d-i)$ for any $0 \le i\le r-1$.  

From the inequalities in Eqs.\,(\ref{ineq-1}) and  
(\ref{ineq-2}), we see that  
Eq.\,(\ref{eqn}) can not hold, which means that our assumption 
$d\geq r$ can not hold. Therefore, we have $o(\Phi)>\deg f(y)$. 
Then $\Phi(\p_y)f(y) = 0$ and, by Lemma \ref{lemma2}, 
we also have $P(x,y) = f(y)$.
\epfv

Now we are ready to prove Theorem \ref{1+P-Case}. \\

\underline{\it Proof of Theorem \ref{1+P-Case}}:\hspace{1ex}
Let $\Lambda = \p_x-\Phi(\p_y)$ and $0\ne P(x,y)\in \bC[x,y]$ such that $\Lambda ^m(P^m)=0$ for any $m\geq 1$. 
Fix any $g(x,y)\in \bC[x,y]$, we want to show that 
$\Lambda ^m(P^mg)=0$ when $m\gg0$.

First, by Lemma \ref{lemma2}, we know $P=e^{x\Phi(\p_y)}(f(y))$ for some $f(y)\in \bC[y]$. Furthermore, from Lemma \ref{lemma1} it follows that either $\Phi(\xi)=0$ or $o(\Phi(\xi))\geq 1$.
If $\Phi=0$, then $\Lambda =\p_x$ and $P=f(y)$.  
In this case it immediately follows that $\Lambda^m(P^mg)=0$ for any $m\geq deg_x\, g$.

So from now on assume that $\Phi\not=0$ and $o(\Phi)\geq 1$. Let $k=deg \, \Phi$. By Eq.\,(\ref{Def-Phi}), we have 
$\Phi(\p_y)=q_1\p_y+ \cdots +q_k \p_y^k$.

If $q_1\not=0$, then we can perform the coordinate change 
$(x,y)\rightarrow (x',y')=(x,y+q_1x)$ to get 
$\tilde{\Lambda}=\p_{x'}-(q_2\p_{y'}^2+\cdots+q_k\p_{y'}^k)$, 
$P(x,y)=P(x',y'-q_1x')=\tilde{P}(x',y')$ and 
$g(x,y) = g(x',y'-q_1x')=\tilde{g}(x',y')$.
Let $\tilde{\Phi}(\p_{y'})=q_2\p_{y'}^2+\cdots+q_t\p_{y'}^t$, then $\tilde{\Lambda}=\p_{x'}-\tilde{\Phi}(\p_{y'})$.

Because that, $\tilde{\Lambda}^m(\tilde{P}^m)=0 \Leftrightarrow \Lambda^m(P^m)=0$ and $\tilde{\Lambda}^m(\tilde{P}^m\tilde{g})=0 \Leftrightarrow \Lambda^m(P^mg)=0$, so we may also assume that $o(\Phi)\geq 2$. Then, by Lemma \ref{lemma3}, we have 
$o(\Phi)>\deg f$ and $P(x,y) = f(y)$. 

Now, for any $m\ge 1$, we have
\allowdisplaybreaks{
\begin{align}
\Lambda^m(f^m(y)g(x,y))&=(\p_x-\Phi(\p_y))^m (f^m(y)g(x, y)) \nno \\
&=\sum_{\substack{a+b=m\\ a, b\ge 0}}(-1)^b \binom mb \p_x^a \Phi^b(\p_y) 
\left(f^m(y)g(x, y)\right) \nno \\
&=\sum_{\substack{a+b=m\\ a, b\ge 0}} (-1)^b \binom mb \Phi^b(\p_y) 
\left(f^m(y)\p_x^a g(x, y)\right). \nno
\end{align} }

Note that, the general term $(-1)^b \binom mb \Phi^b(\p_y) 
\left(f^m(y)\p_x^a g(x, y)\right)$ in the summation above 
is not equal to zero only if 
\begin{align}
\begin{cases}
a\le \deg_x g,\\
b\, o(\Phi)=(m-a) o(\Phi) \le m\deg f+\deg_y g.
\end{cases}
\end{align}

Assume the inequalities above hold, Then, 
we have 
\begin{align}
m(o(\Phi)-\deg f)\le a \, o(\Phi)+\deg_y g 
\le o(\Phi) \deg_x g+\deg_y g.
\end{align}

Since $o(\Phi)>\deg f$ (as pointed out above), 
the combined inequality above is equivalent to 
\begin{align}
m\le  \frac{o(\Phi) \deg_x g+\deg_y g}{o(\Phi)-\deg f}.
\end{align}
Therefore, from the arguments above, we have 
$\Lambda^m(f^mg)=0$ for each  
\begin{align}
m>\frac{o(\Phi) \deg_x g+\deg_y g}{o(\Phi)-\deg f}.
\end{align}
Hence, we have proved Theorem \ref{1+P-Case}.
\epfv

Next we give the following example to show that, 
for the two variable case, both Conjectures \ref{VC} 
and \ref{GVC} fail when $P(x, y)$ is allowed to be 
a formal power series, instead of just a polynomial.

\begin{exam}\label{Counter-GVC}
Let $\Lambda=\p_y \p_x$ and $P(x, y)=x+e^y$. 
Then, for any $m\ge 1$, we have
\begin{align*}
\Lambda^m(P^m)&=\p_y^m \p_x^m (x+e^y)^m =\p_y^m (m!)=0,\\
\Lambda^m(P^{m+1})&=\p_y^m \p_x^m (x+e^y)^{m+1} =(m+1)! 
\, \p_y^m (x+e^y)\\
&=(m+1)!\, e^y\ne 0.
\end{align*}
Hence Conjecture \ref{VC} fails in this case.

Furthermore, let $g(x, y)=x$. Then, for any $m\ge 1$,  we have 
\begin{align*}
\Lambda^m(P^mg)&=\p_y^m \p_x^m (x(x+e^y)^m) 
=\p_y^m\left(m!x + \binom m1 \p^{m-1} (x+e^y)^m\right) \\
&=m!\, \p_y^m(x + m(x+e^y))=m\,m!\, e^y\ne 0.
\end{align*}
Hence, Conjecture \ref{GVC} also fails for the formal power series $P(x, y)=x+e^y$.
\end{exam}

\renewcommand{\theequation}{\thesection.\arabic{equation}}
\renewcommand{\therema}{\thesection.\arabic{rema}}
\setcounter{equation}{0}
\setcounter{rema}{0}

\section{\bf Some Results on Rational Polytopes and the Density Theorem of Polytopes of Laurent Polynomials}\label{S4} 

In this section, we first recall and prove some results on rational polytopes of $\bR^n$ that will be needed later in this paper. We then use 
the Duistermaat-van der Kallen theorem, 
Theorem \ref{ThmDK}, to prove   
what we call the {\it density theorem} of polytopes of Laurent polynomials (See Theorem \ref{DensityThm}) along with   
some of its variations.  
We will also show in Lemma \ref{Pnot>L} that Conjecture 
\ref{GVC} holds when the polytope 
$\poly(P)-\poly(\Lambda)$ has 
no intersection points with 
$(\bR^{\ge 0})^{\times n}$. Finally, 
we give an example, Example \ref{Counter-ThmDK}, 
to show that the Duistermaat-van der Kallen theorem 
can not be generalized to the setting of 
Laurent formal power series.

First, let us fix the following notations and conventions which 
together with the notations fixed in Section \ref{S1}  
will be used throughout the rest of this paper.  \\

{\bf Notation and Convention:} \\

\begin{enumerate}
\item We use $x=(x_1, x_2, ..., x_n)$ to denote the coordinates of the Euclidean space $\bR^n$. For any $u=(a_1, a_2, ..., a_n) \in \bR^n$,  
against the traditional notation we set 
$|u|\!:=\sum_{i=1}^n a_i.$ Note that $|u|$ 
could be negative for some $u\in \bR^n$.

\item For any non-zero $u\in \bR^n$, we denote by $R_u$ 
the ray with the $0\in \bR^n$ as its (only) end point and 
passing through $u$. When $u=0$, we let $R_u$ denote 
the single point $0\in \bR^n$.  

\item For any $m\in \bZ$, $\beta\in \bR^n$ and a subset 
$A\subset \bR^n$, we set
\begin{align*}
mA\!:=&\{ mu\,|\, u \in A\},\\
\beta\pm A\!:=&\{\beta\pm u\,|\, u\in A \}.
\end{align*}

\item We denote by $\qrt$ the set of the 
vectors in $\bR^n$ whose components are 
all non-negative. Furthermore, 
we introduce a partial order $>$ for 
vectors in $\bR^n$ by setting, for any $u, v\in \bR^n$, 
$u>v$ if $u\ne v$ and $u-v\in \qrt $; and 
$u\ge v$ if $u>v$ or $u=v$.

\item For any finite subset 
 $A=\{u_i\,|\, 1\le i\le k\} \subset \bR^n$, we set   
\begin{align}\label{polytope}
\poly(A)\!:=\left\{\sum_{i=1}^k c_i u_i 
\,\left.\right|\, c_i\ge 0;\, \sum_{i=1}^k c_i=1 \right\}.
\end{align}
We call the subset above the {\it polytope} generated by 
the points $u_i\in \bR^n$ $(1\le i\le k)$, or simply, 
by the subset $A$. Throughout the paper, 
by a {\it polytope} we always mean a subset of $\bR^n$ of 
the form as in Eq.\,$(\ref{polytope})$.

\item For any $u\in \bR^n$, we say $u$ is a {\it rational point} of $\bR^n$ if its all components are rational numbers. We say a polytope is {\it rational} if all its vertices are rational.

\item For any fixed Laurent polynomial $P(z)\in \czz$, 
and any $\alpha\in\bZ^n$, we denote by 
$[z^\alpha]P(z)$ the coefficient of $z^\alpha$ in $P(z)$. 
We define the {\it support} of $P(z)$, denoted by 
$\supp (P)$, to be the subset of $\alpha \in \bZ^n$ 
such that $[z^\alpha]P(z)\ne 0$; and 
the {\it polytope} of $P(z)$, denote by $\poly (P)$, 
to be the polytope generated by $\supp(P)$.

\item For any differential operator 
$\Lambda=\Lambda(\p)$ with $\Lambda(\xi)\in \bC[\xi]$, 
we define the {\it support} of $\Lambda$,
denoted by $\supp(\Lambda)$, and the {\it polytope} 
of $\Lambda$, 
denoted by $\poly(\Lambda)$, to be the support and 
the polytope, respectively, of the 
polynomial $\Lambda(\xi)$.
\end{enumerate}

We start with the following lemma which is well-known
(e.g. see \cite{CLO}) and also easy to prove directly. 
In our later argument we will frequently use this lemma 
without explicitly referring to it. 

\begin{lemma}\label{NewlyAdded}
$(a)$ For any polytope $\Sigma$ in $\bR^n$ 
and $m\ge 1$, we have
\begin{align}
m\Sigma=\left\{ \sum_{i=1}^m u_i \,|\, u_i\in \Sigma\right\}.
\end{align} 
$(b)$ For any Laurent polynomial $P(z)\in \czz$  
and $m\ge 1$, we have
\begin{align}
\poly(P^m)=m\poly(P).
\end{align} 
\end{lemma}

The following lemma is also well-known. 
But, for the sake of completeness, 
we include a proof here.
 
\begin{lemma}\label{MoveAway}
For any $\beta\in\bR^n$ and any 
polytope $\Sigma$ in $\bR^n$ with  
$\Sigma\cap\qrt = \emptyset$, 
there exists $N\ge 1$ such that
$(\beta+m\Sigma)\cap \qrt =\emptyset$ 
for any $m\ge N$.
\end{lemma}
\pf 
Assume otherwise, then there exist a strictly increasing 
sequence $\{m_k \,|\, m_k\in \bN\}$ and a sequence 
$\{v_k \, |\, v_k\in (\beta+m_k\Sigma)\cap \qrt \}$.

Note that, for any $k\ge 1$, we may write 
$v_k=\beta+m_k u_k$ for some $u_k\in \Sigma$.
Since $\Sigma$ is bounded and closed, and hence compact,
the sequence $\{u_k\,|\, k\ge 1\}$ has a subsequence which 
converges to an element $u\in \Sigma$. 
Without losing any generality,  
we still denote this subsequence by 
$\{u_k\,|\, k\ge 1\}$.  

Note that, for any $k\ge 1$, $u_k=v_k/m_k-\beta/m_k$. 
Then we have  
\begin{align*}
u=\lim_{k\to \infty} u_k= \lim_{k\to \infty} 
(v_k/m_k-\beta/m_k)= \lim_{k\to \infty} (v_k/m_k).
\end{align*}

Since, for any $k\ge 1$,  $m_k\ge 1$ and $v_k\in \qrt $, 
we have $v_k/m_k\in \qrt$. Furthermore, 
since $\qrt$ is closed, from the equation above 
we have $u\in \qrt $. Therefore, we have 
$u\in  \Sigma \cap \qrt$, which is a 
contradiction. Hence the lemma holds. 
\epfv

\begin{lemma}\label{QD-polytope}
For any two rational polytopes $\Sigma$ and $\Gamma$, we have  

$(a)$ If $\Sigma \cap \Gamma\ne \emptyset$, then it is 
also a rational polytope.

$(b)$ $\Sigma-\Gamma\!:=\{ u-v\,|\, u\in \Sigma,\, v\in \Gamma\}$ 
is also a rational polytope.

$(c)$ For any rational $w\in \Sigma-\Gamma$, there exist 
$u\in \Sigma$ and $v\in \Gamma$ such that $u$, $v$ are both  rational and $w=u-v$. 
\end{lemma}

\pf $(a)$ It is well known 
(e.g. see Theorem 1.1, pp.\,29 in \cite{Zi}) 
that any polytope is 
a set of common solutions of a system of 
linear equations or inequalities. 
It is easy to see that  
a polytope is rational iff its determining 
linear equations or inequalities are defined over 
$\bQ$, i.e. all the coefficients of unknowns including 
constant terms of the equations or 
inequalities are in $\bQ$.  
Since $\Sigma \cap \Gamma$ is determined by the union 
of determining linear equations or inequalities
of $\Sigma$ and $\Gamma$, 
$\Sigma \cap \Gamma$ is also a rational polytope   
if it is not empty.

$(b)$ First, it is easy to see that 
$(-1)\Gamma$ is a rational polytope and 
$\Sigma-\Gamma$ is the same as the so-called 
Minkowski sum of the polytopes $\Sigma$ 
and $\Gamma$. It is well known that 
the Minkowski sum of any two polytopes is also 
a polytope (e.g. see \cite{CLO} and \cite{Zi}).  
The reason that the polytope $\Sigma-\Gamma$ is 
also rational is because 
 any vertex of $\Sigma-\Gamma$ is 
the difference of a vertex of $\Sigma$ 
and a vertex of $\Gamma$.

$(c)$ Note that, for the fixed 
$w\in \Sigma-\Gamma$ in the lemma, 
the set of elements $u\in \Sigma$ such that 
$w=u-v$ for some $v\in \Gamma$ is given 
by $\Sigma \cap (w+\Gamma)$ 
which is non-empty by the existence 
of $w$ itself.

Since both $w$ and $\Gamma$ are 
rational, so is $w+\Gamma$. 
By $(a)$ we know that 
$\Sigma \cap (w+\Gamma)$ is 
a (non-empty) rational polytope 
since $\Sigma$ is also rational.

Let $u$ be any rational point of 
$\Sigma \cap (w+\Gamma)$, say a vertex of 
this polytope, and write it as $u=w+v$ for 
some $v\in \Gamma$. Since $v=u-w$, 
$v$ is also rational. 
Hence we get $(c)$.
\epfv

\begin{corol}\label{Q-polytope}
For any rational polytope $\Sigma$ of $\bR^n$ with 
$\Sigma\cap\qrt\ne \emptyset$, there exists a rational 
$u\in \Sigma\cap\qrt$.
\end{corol}

\pf Since the rational polytope $\Sigma$ 
is a closed and bounded subset of $\bR^n$, so is 
$\Sigma\cap\qrt$. We may choose an $N\in \bN$ 
such that $\Sigma\cap\qrt$ lies inside 
the polytope generated by $0\in \bR^n$ and 
$Ne_i$ $(\lin)$, where $e_i$'s are 
the vectors in the standard basis 
of $\bR^n$. Since the latter polytope 
is also rational, the corollary follows 
from Lemma \ref{QD-polytope}, $(a)$. 
\epfv

Next we prove some results for the polytopes 
of Laurent polynomials. But, first let us recall 
the following remarkable theorem which was first 
conjectured by O. Mathieu \cite{Ma} and later 
was proved by J. Duistermaat and 
W. van der Kallen \cite{DK}. 

\begin{theo}\label{ThmDK} $(${\bf  Duistermaat and van der Kallen}$)$ 
For any $f(z)\in \czz$ such that the constant term of 
$f^m(z)$ is equal to zero for any $m\ge 1$, we have   
$0\not \in \poly(f)$.
\end{theo}

The following result will play some crucial roles in later sections. We believe that it is also important in its own right, so we formulate it as a theorem and call it the {\it density theorem} of polytopes of Laurent polynomials.

\begin{theo}\label{DensityThm}$(${\bf The Density Theorem}$)$ 
For any $P(z)\in \czz$ and any rational point $u\in \poly(P)$, 
there exists $m\ge 1$ such that 
$R_u\cap \supp(P^m)\ne \emptyset$. 
\end{theo}

\pf Assume otherwise, i.e. 
$R_u\cap \supp(P^m)=\emptyset$ for any $m\ge 1$.

Since $u$ is rational, there exists $N\ge1$ 
such that $\beta:=N u\in \bZ^n$. Hence, we have 
$R_u=R_\beta$. 
Since $u\in \poly(P)$ and $\poly(P^N)=N\poly(P)$, 
we have 
\begin{align}\label{DensityThm-pe1}
\beta=Nu\in\poly(P^N).  
\end{align}

Let $f(z)\!:=z^{-\beta} P^N(z)\in \czz$. Then, 
for any $m\ge 1$, the constant term of $f^m$ 
is equal to zero. 
Otherwise, we would have $m\beta\in \supp(P^{mN})$ 
and $m\beta \in R_\beta \cap \supp(P^{mN})
=R_u\cap \supp(P^{mN})$, which contradicts 
to our assumption on $R_u$.

Now apply Theorem \ref{ThmDK} to $f(z)$, we have  
$0\not \in \poly(f)$. But, on the other hand, we have  
\begin{align*}
\poly(f)=-\beta+\poly(P^N).
\end{align*}
Therefore, $\beta\not \in \poly(P^N)$
which contradicts to Eq.(\ref{DensityThm-pe1}).
Hence, we have proved the theorem.
\epfv

From the proof above, we can actually get a 
stronger result.

\begin{corol}\label{DensityThm-3}
For any $P(z)\in \czz$ and any rational point 
$u\in \poly(P)$, there exist infinitely many 
$m_i \ge 1$ such that 
\begin{align*}
R_u\cap \supp(P^{m_i})\ne \emptyset.
\end{align*} 
\end{corol}

\pf Assume otherwise, then there exists 
$N\ge 1$ such that $R_u \cap \supp(P^{m})
= \emptyset$ for any $m\ge N$. 

But, applying  
Theorem \ref{DensityThm} to the Laurent polynomial
$P^N(z)$ and the rational point 
$Nu\in \poly(P^N)$, we see that there exists $m\ge 1$, 
such that 
\begin{align*}
\emptyset\neq R_{Nu}\cap \supp(P^{mN})
= R_{u}\cap \supp(P^{mN}).
\end{align*}
Since $mN\ge N$, we get a contradiction.
\epfv

When $P(z)$ is homogeneous with respect to
the generalized degree (counting 
$\deg z_i^{-1}=-1$ for any $\lin$) 
of Laurent polynomials, 
we have the following more 
precise result.

\begin{corol}\label{DensityThm-2}
For any homogeneous Laurent polynomial $P(z)$ 
of degree $d\ne 0$, and any rational point $u\in \poly(P)$, there exist infinitely many 
$m_i \ge 1$ such that  $m_iu\in \supp(P^{m_i})$. 
\end{corol}
\pf First, since $\deg P=d$, we see that 
$\poly(P)$ lies in the affine hyperplane $H$ determined 
by the equation $\sum_{i=1}^n x_i=d$. Since 
$u\in \poly(P)\subset H$, we have $|u|=d$.

Second, by Corollary \ref{DensityThm-3}, we know that 
there exist infinitely many $m_i \ge 1$ such that 
$R_u\cap \supp(P^{m_i})\ne \emptyset$. 
We fix any such an $m_i\ge 1$ and choose any 
$\beta_i\in R_u\cap \supp(P^{m_i})$ (Actually, such $\beta_i$ is unique as the intersection point of $R_u$ with the affine hyperplane $m_iH$). Write $\beta_i=k_iu$ for some $k_i\ge 0$. 
Then we have $|\beta_i|=k_i|u|=k_id$. Since 
$\beta_i\in \supp(P^{m_i})\subset m_i\poly(P)\subset m_iH$,
we have $|\beta_i|=m_id$. Hence we have $k_i=m_i$ and 
$\beta_i=m_iu$.
\epfv

Next we prove the following general 
result on Conjecture \ref{GVC}, which will be needed in Section \ref{S6}.

\begin{lemma}\label{Pnot>L}
Let $P(z)\in \cz$ and $\Lambda=\Lambda(\p)$ 
any differential operator with constant coefficients. 
Assume
\begin{align}\label{Pnot>L-e1}
(\poly(P)-\poly(\Lambda))\cap \qrt=\emptyset.
\end{align}
Then, we have

$(a)$ for any $m\ge 1$, $\Lambda^m (P^m)=0$.

$(b)$ for any $g(z)\in \cz$,  
$\Lambda^m (P^mg)=0$ when $m\gg0$.
\end{lemma}

\pf We first prove $(b)$ as follows.

First, note that, by the linearity of 
$\Lambda^m (P^mg)$ $(m\ge 1)$ on $g(z)\in \cz$, 
it is easy to see that we may assume 
$g(z)=z^\gamma$ for some 
$\gamma\in \bN^n$.

Second, let $\Sigma\!:=\poly(P)-\poly(\Lambda)$. 
Then, we have that $\Sigma\cap\qrt=\emptyset$ and also, 
by Lemma \ref{QD-polytope}, $(b)$, $\Sigma$ is a polytope 
of $\bR^n$. 

Apply Lemma \ref{MoveAway} to the polytope 
$\Sigma$ and $\gamma\in \bN^n$, we get a 
$N\ge 1$ such that, for any $m \ge N$,
\begin{align}\label{Pnot>L-pe1}
(\gamma+m\Sigma)\cap \qrt=\emptyset.
\end{align}

On the other hand, for any $m\ge 1$, we have
\begin{align}
m\Sigma=m\poly(P)-m\poly(\Lambda)=\poly(P^m)-\poly(\Lambda^m).
\end{align} 
Therefore, by Eq.(\ref{Pnot>L-pe1}) and the equation above,
we have 
\begin{align}
\qrt\cap \left(\gamma +(\poly(P^m)-\poly(\Lambda^m) \right)=\emptyset
\end{align}
for any $m\ge N$.

Consequently, for any $m\ge N$, 
$\lambda \in \supp(P^m)\subset \poly(P^m)$ and 
$\mu \in \supp(\Lambda^m) \subset \poly(\Lambda^m)$, 
we have 
\begin{align}
(\gamma+\lambda)-\mu  
=\gamma+(\lambda-\mu)\not \in \qrt.
\end{align} 
Hence, we have $\p^\mu z^{\gamma+\lambda}=0$.
Since $\Lambda^m$ is a linear 
combination of 
$\p^\mu$ $(\mu\in \supp(\Lambda^m))$ and 
$z^\gamma P^m(z)$ is a linear combination 
over $\bC$ of $z^{\gamma+\lambda}$ 
$(\lambda\in \supp(P^m))$, 
we have $\Lambda^m(z^\gamma P^m(z))=0$ for 
any $m\ge N$. Therefore, we have proved $(b)$.

To see $(a)$ also holds, note that, by choosing 
$\gamma=0\in\bN^n$ in the proof above,  
Eq.(\ref{Pnot>L-pe1}) actually holds for any $m\ge 1$.
This is because the condition 
$\Sigma\cap\qrt=\emptyset$ implies directly 
$m\Sigma\cap\qrt=\emptyset$ for any $m\ge 1$. 
Therefore, the argument above also goes 
through with $N=1$ and $\gamma=0$, 
which means $(a)$ also holds.   
\epfv

\begin{rmk}\label{Rrk4.1.1}
Note that, by Lemma \ref{Pnot>L}, $(a)$, 
we see that the condition given in 
Eq.\,$(\ref{Pnot>L-e1})$ implies the condition in  
Conjecture \ref{GVC}. By Lemma \ref{Pnot>L}, $(b)$, we see  
that Conjecture \ref{GVC} does hold in this case.
\end{rmk}

Finally, let us point out that, like Conjectures \ref{VC} 
and \ref{GVC}, the Duistermaat-van der Kallen theorem, 
Theorem \ref{ThmDK} can not be 
generalized to Laurent formal 
power series either.

\begin{exam}\label{Counter-ThmDK}
Let $f(x, y)=y^{-1}(1+x^{-1}e^y)$ and $g(x, y)=x$. Then, for any $m\ge 1$, it is easy to check that, the constant term of $f^m$ is equal to zero, but the constant term of 
$f^mg$ is equal to $1/(m-1)!\ne 0$. Therefore, 
$0\in \supp(f^mg)$ for each $m\ge 1$.

Assume that the statement of 
Theorem \ref{ThmDK} holds for $f(x, y)$, i.e.
$0\not \in \poly(f)$, then, as shown in \cite{DK}, we will have $0\not \in \supp(f^mg)$ when $m\gg0$. But this contradicts to the fact pointed out above, i.e.\@ 
$0\in \supp(f^mg)$ for each $m\ge 1$. 
Hence the statement of Theorem \ref{ThmDK} 
fails for the Laurent power series $f(x, y)$.
\end{exam} 

\renewcommand{\theequation}{\thesection.\arabic{equation}}
\renewcommand{\therema}{\thesection.\arabic{rema}}
\setcounter{equation}{0}
\setcounter{rema}{0}

\section{\bf A Theorem on Laurent Polynomials with no Holomorphic Parts 
and the Vanishing Conjecture when $\Lambda(\p)$ or $P(z)$  
is a Monomial}\label{S5}

In this section, we first prove a conjecture proposed 
in \cite{GVC} on Laurent polynomials with 
no holomorphic parts (See Theorem \ref{H-Case}). 
By using this result, we then show in 
Corollary \ref{M-Case}
that Conjecture \ref{GVC} 
holds when the polynomials $P(z)$ or the 
differential operators $\Lambda=\Lambda(\p)$ 
is a monomial of $z$ or $\p$, respectively.

\begin{theo}\label{H-Case}
Let $f(z)\in \czz$ such that, for any $m\ge 1$, 
the holomorphic part of $f^m$ is equal to zero, 
i.e. $\supp(f^m)\cap\qrt=\emptyset$, or equivalently, 
$[z^\alpha]f^m=0$ for any $\alpha\in \bN^n$.
Then, for any $g(z)\in \czz$, 
the holomorphic part of $g f^m$ is equal to zero 
when $m\gg0$.
\end{theo}

\pf Since the subset of 
Laurent polynomials with no holomorphic 
parts is a subspace of $\czz$,  
and $g(z)$ is a linear combination over $\bC$ 
of finitely many monomials of $z$, it is easy to 
see that we may assume $g(z)=z^\beta$ 
for some $\beta\in \bZ^n$.  

Let us first show $\poly(f)\cap\qrt =\emptyset$ 
by the contradiction method. 

Assume otherwise, i.e. $\poly(f)\cap\qrt \ne \emptyset$. 
Since $\poly(f)$ is a rational polytope, by Corollary \ref{Q-polytope}, there exists a rational $u\in \poly(f)\cap\qrt$. 
Below we fix any such a rational point $u$.

Apply Theorem \ref{DensityThm} 
to $f(z)\in \czz$ and the rational point 
$u\in \poly(f)$, there exists $N\ge 1$ such that 
$R_u\cap \supp(f^N)\ne \emptyset$.
Note that $R_u\subset \qrt$ since 
$u\in \qrt$. So we have   
$\qrt \cap \supp(f^N)\ne \emptyset$,   
and hence the holomorphic part of $f^N$ 
is not zero, which is a contradiction.
Therefore, we must have 
$\poly(f)\cap\qrt =\emptyset$. 

Next, apply Lemma \ref{MoveAway} to the polytope 
$\poly(f)$ and $\beta\in \bN^n$,   
we know that 
$\qrt \cap (\beta+m\poly(f))=\emptyset$ 
when $m\gg0$. Since, for any $m\ge 1$, 
 $m\poly(f)=\poly(f^m)$ and 
$\supp(f^m)\subset \poly(f^m)$,  
we have $\qrt \cap (\beta+\supp(f^m))=\emptyset$ 
when $m\gg0$. 

On the other hand, it is easy to see that 
$\supp(z^\beta f^m)=\beta+\supp(f^m)$ for 
any $m\ge 1$. Therefore, we have 
$\qrt \cap \supp(z^\beta f^m)=\emptyset$ 
when $m\gg0$, which means that the 
theorem holds for $g(z)=z^\beta$. 
\epfv

Next, we use Theorem \ref{H-Case} to 
show that Conjecture \ref{GVC} holds when 
$P(z)$ or $\Lambda$ is a monomial of 
$z$ or $\p$, respectively. But, first, let us formulate 
the following simple observation as a lemma since it
will be crucial for our later arguments.

\begin{lemma}\label{TrivialLemma}
For any $\alpha, \beta\in \bN^n$, we have that, 
$\p^\alpha z^\beta=0$ iff $\beta\not \ge \alpha$, 
or equivalently, $\beta-\alpha \not \in \qrt$.
\end{lemma}

Note that the lemma above is not necessarily 
true for $\beta\in\bZ^n\backslash\bN^n$. 

\begin{corol}\label{M-Case}
Conjecture \ref{GVC} holds if either $\Lambda=\p^\alpha$ 
or $P(z)=z^\alpha$ for some 
$\alpha\in \bN^n$. 
\end{corol}

\pf We prove the corollary for the case that 
$P(z)=z^\alpha$ for some $\alpha\in \bN^n$. 
The proof for the other case is similar.

First, by the linearity on $g(z)\in \cz$, 
we may assume that $g(z)=z^\gamma$ 
for some $\gamma\in \bN^n$.

Second, for any $h(\xi)\in \bC[\xi]$ and any $\beta\in \bN^n$, 
by Lemma \ref{TrivialLemma}, 
it is easy to see that we have the following equivalences:
\begin{align}\label{M-Case-pe1}
h(\p) z^{\beta}=0 \quad &\Leftrightarrow  \quad 
\beta \not \ge \mu \mbox{ for any } \mu\in \supp(h)   \\
&\Leftrightarrow  \quad 
\supp(z^\beta h(z^{-1})) \cap\qrt=\emptyset.\nno
\end{align}
Note that the last statement above is equivalent to saying 
that the holomorphic part of the Laurent polynomial
$z^\beta h(z^{-1})$ is equal to zero. 

Now, we write $\Lambda=\Lambda(\p)$ for some 
polynomial $\Lambda(\xi)\in \bC[\xi]$ and set 
$f(z)\!:=\Lambda(z^{-1})z^\alpha$. 
Since $\Lambda^m (P^m)=\Lambda^m (z^{m\alpha})=0$ 
for any $m\ge 1$, applying the equivalences 
in Eq.\,(\ref{M-Case-pe1}) 
with $h(\xi)=\Lambda^m(\xi)$ and 
$\beta=m\alpha$, we know that  
the holomorphic part of 
$f^m=\Lambda^m(z^{-1})z^{m\alpha}$ is equal to
zero for any $m\ge 1$. 

Applying Theorem \ref{H-Case} to $f(z)$ and 
$g(z)=z^\gamma$, we know that there exists a $N\ge 1$ 
such that the holomorphic part of 
$f^m z^\gamma=\Lambda^m(z^{-1})z^{m\alpha+\gamma}$ is 
equal to zero for any $m\ge N$.
For any fixed $m\ge N$, 
applying the equivalences 
in Eq.\,(\ref{M-Case-pe1}) 
with $h(\xi)=\Lambda^m(\xi)$ and 
$\beta=m\alpha+\gamma$, we get 
$\Lambda^m (P^m z^\gamma)
=\Lambda^m (z^{m\alpha+\gamma})=0$.
Therefore, Conjecture \ref{GVC} does 
hold in this case.  
\epfv

\renewcommand{\theequation}{\thesection.\arabic{equation}}
\renewcommand{\therema}{\thesection.\arabic{rema}}
\setcounter{equation}{0}
\setcounter{rema}{0}

\section{\bf Proof of the Vanishing Conjecture for the Differential Operator 
$\Lambda=a\p^\alpha+b\p^\beta$ with 
$|\alpha|\ne |\beta|$}\label{S6}

In this section, we prove the following case of 
Conjecture \ref{GVC}.

\begin{theo}\label{MM-Case}
Let $\alpha, \beta \in \bN^n$ such that 
$|\alpha|\ne |\beta|$, and $\Lambda=a\p^\alpha+b\p^\beta$
for some $a, b\in\bC$. Then, Conjecture \ref{GVC} holds for the differential operator $\Lambda$ and 
any homogeneous $P(z)\in \cz$.
\end{theo}

First, let us consider the following simple cases of the theorem above.

If $a=b=0$, then $\Lambda=0$. There is nothing to prove. 
If only one of the $a$ and $b$ is zero, then, after a change 
of variables, we may assume 
$\Lambda=\p^\gamma$ for some $\gamma\in \bN^n$. 
Then Conjecture \ref{GVC} in this case follows 
directly from Corollary \ref{M-Case}. 

Therefore we may assume that $a$ and $b$ are both nonzero.
Then, by using the fact that $\alpha\ne \beta$, 
it is easy to see that, after a change of variables, 
we may assume $a=b\ne 0$. Note also that in general    
Conjecture \ref{GVC} holds for a  
differential operator iff it holds for any  
nonzero scalar multiple of the differential operator. 
So we may further assume that $a=b=1$. 
Therefore, we can reduce Theorem \ref{MM-Case} 
to the case when $\Lambda=\p^\alpha+\p^\beta$ for some 
$\alpha, \beta\in \bN^n$ with 
$|\alpha| \ne |\beta|$. 

Throughout the rest of this section, 
we will fix a differential operator $\Lambda$ 
as above and a homogeneous polynomial 
$P(z)\in \cz$ of degree $d=\deg P\ge 0$ 
such that $\Lambda^m (P^m)=0$ for any $m\ge 1$. 
We divide the proof of Theorem \ref{MM-Case} 
for this case into several lemmas. 

\begin{lemma}\label{lemma-MM1}
$(a)$ For any $m\ge 1$, we have
\begin{align}
\supp(\Lambda^m)=\{k\alpha+\ell \beta\,|\, k, \ell \in \bN; 
k+\ell=m\}.
\end{align}

$(b)$ For any $m, k\ge 1$ and $u\in \supp(\Lambda^k)$, 
 $mu\in \supp(\Lambda^{mk})$.
\end{lemma}

\pf $(a)$ The statement follows directly from the binomial expansion of $\Lambda^m=(\p^\alpha+\p^\beta)^m$. 

$(b)$ By $(a)$, we may write 
$u=r\alpha+s\beta$ for some $r, s\in \bN$ 
with $r+s=k$. Then, we have 
$mu= mr\alpha+ms\beta$. Since $mr, ms\in \bN$ 
and $mr+ms=mk$, by $(a)$ again we have 
$mu\in \supp(\Lambda^{mk})$.
\epfv

\begin{lemma}\label{lemma-MM2}
For any fixed $m\ge 1$, we have

$(a)$ $\p^\mu P^m=0$ for any $\mu\in \supp(\Lambda^m)$.

$(b)$ for any $\mu\in \supp(\Lambda^m)$ and
$\gamma \in \supp(P^m)$, we have 
$\gamma\not \ge \mu$ or equivalently,
$\gamma-\mu \not \in \qrt$.
\end{lemma}

\pf $(a)$ We first consider
\begin{align}\label{lemma-MM1-pe1}
0&=\Lambda^m (P^m)=(\p^\alpha+\p^\beta)^m (P^m)
=\sum_{\substack{k, \ell \ge 0\\k+\ell=m}}\binom mk 
\p^{k\alpha+\ell\beta} (P^m).
\end{align}
Note that, for any $k, \ell \ge 0$ with 
$k+\ell=m$, we have
\begin{align}\label{lemma-MM1-pe2}
\deg(\p^{k\alpha+\ell\beta} P^m)=
md-(k|\alpha|+\ell|\beta|).
\end{align}
Assume that, for some $k', \ell' \in \bN^n$ with 
$k'+\ell'=m$, we have 
\begin{align*}
\deg(\p^{k' \alpha+\ell' \beta} P^m)
=\deg(\p^{k\alpha+\ell\beta} P^m)
\end{align*}
Then, by Eq.\,(\ref{lemma-MM1-pe2}), we get  
\begin{align}
k|\alpha|+\ell|\beta|&=k'|\alpha|+\ell'|\beta|, \nno \\
(k-k')|\alpha|&=(\ell'-\ell)|\beta|. \label{lemma-MM1-pe3}
\end{align}

Since $k+\ell=k'+\ell'=m$, we have $k-k'=\ell'-\ell$.
Combining Eq.\,(\ref{lemma-MM1-pe3}) with the fact that 
$|\alpha|\ne |\beta|$, 
we see that $k-k'=\ell'-\ell=0$. 
Therefore, all the terms in the sum of 
Eq.\,(\ref{lemma-MM1-pe1}) have different degrees. 
Hence they all have to be zero.

On the other hand, by Lemma \ref{lemma-MM1}, $(a)$, any 
$\mu\in \supp(\Lambda^m)$ has the form 
$k\alpha+\ell \beta$ for some $k, l\in \bN$ 
with $k+l=m$. Hence $(a)$ follows.

$(b)$ First, we write $P^m(z)$ as 
$P^m(z)=\sum_{\gamma\in \supp(P^m)} b_\gamma z^\gamma$
with $b_\gamma\in \bC^{\times}$. 
For any $\mu\in \supp(\Lambda^m)$, 
by $(a)$ and also Lemma \ref{TrivialLemma}, 
we have
\begin{align}\label{lemma-MM1-pe4}
0=\p^{\mu}(P^m) = 
\sum_{\gamma\in \supp(P^m)} b_\gamma 
\p^\mu (z^\gamma)  
= \sum_{\substack{\gamma\in \supp(P^m)\\
\gamma \ge \mu }} 
 b_{\gamma} z^{\gamma-\mu}.
\end{align}

Note that $b_\gamma\ne 0$ for any $\gamma\in \supp(P^m)$ and, 
for any $\gamma_1\ne \gamma_2 \in \supp(P^m)$, 
$z^{\gamma_1-\mu}\ne z^{\gamma_2-\mu}$. Then, 
from Eq.\,(\ref{lemma-MM1-pe4}) we see that  
there can not be any $\gamma\in \supp(P^m)$ 
such that $\gamma\ge \mu$. 
Hence $(b)$ also holds. 
\epfv

\begin{lemma}\label{lemma-MM3}
Let $u\in \poly(P)$ and $v\in \poly(\Lambda)$. 
Assume that both $u$ and $v$ are rational. Then, 
we have, $u\not \ge v$.  
\end{lemma}

\pf First, denote by 
$L_{\alpha, \beta}$ the line segment in $\bR^n$ 
connecting $\alpha$ and $\beta$. Then 
it is easy to see that $\poly(\Lambda)=L_{\alpha, \beta}$.
Since both $\alpha$ and 
$\beta$ are rational, 
it is easy to check that,  
a point of $L_{\alpha, \beta}$ is rational 
iff it is a linear combination of $\alpha$ 
and $\beta$ with rational coefficients.
Therefore, we may write $v=r\alpha+s\beta$ 
for some rational $r, s\ge 0$ with 
$r+s=1$. Let $N\ge 1$ such that $Nr, Ns\in \bN$.
Since $Nr+Ns=N$, by Lemma \ref{lemma-MM1}, $(a)$, 
we have $Nv=Nr\alpha+Ns\beta \in \supp(\Lambda^N)$.

Second, apply Corollary \ref{DensityThm-2}
to the homogeneous polynomial $P^N$ and 
the rational point $Nu\in \poly(P^N)$, there exists
$m\ge 1$ such that $mNu\in \supp(P^{mN})$.
Since $Nv\in \supp(\Lambda^N)$, by Lemma \ref{lemma-MM1}, 
$(b)$, we have $mNv\in \supp(\Lambda^{mN})$. 

Note that $mNu\in \supp(P^{mN})$ 
and $mNv\in \supp(\Lambda^{mN})$ are both rational,  
by Lemma \ref{lemma-MM2}, $(b)$, we have 
$mNu-mNv=mN(u-v)\not \in \qrt$. Since $mN>0$, 
we also have $(u-v)\not \in \qrt$, i.e. 
$u\not \ge v$.
\epfv

Now, we can prove the main result, Theorem \ref{MM-Case}, of this section as follows.  

\medskip

\underline{\it Proof of Theorem \ref{MM-Case}}: 
First, by the reductions given at the beginning of this section, 
we may assume that the differential operator 
$\Lambda=\p^\alpha+\p^\beta$ for some $\alpha, \beta\in \bN^n$ 
with $|\alpha| \ne |\beta|$. 

Second, let $\Sigma\!:=\poly(P)-\poly(\Lambda)$. Then, 
by Lemma  \ref{Pnot>L}, it will be enough to show that
$\Sigma \cap \qrt=\emptyset$. 

We assume otherwise, i.e. $\Sigma \cap \qrt\neq \emptyset$, 
and derive a contradiction as follows.

Note first that, by Lemma \ref{QD-polytope}, $(b)$, 
we know that $\Sigma$ is a rational polytope. 
By Corollary \ref{Q-polytope}, there exists a rational point  
$w\in \Sigma\cap\qrt$. Then, by lemma \ref{QD-polytope}, $(c)$, 
there exist $u\in \poly(P)$ and $v\in \poly(\Lambda)$ 
such that $u$, $v$ are both rational and 
$u-v=w \in \Sigma\cap\qrt\subset \qrt$. 
But this contradicts to Lemma \ref{lemma-MM3}. Hence 
we have proved the theorem.
$\Box$ \\

Finally, let us point out that, 
by similar arguments as in the proof of Theorem \ref{MM-Case}, 
it is easy to see that Conjecture \ref{GVC} also holds 
for the following case.  

\begin{corol}\label{MM-Case2}
Let $\alpha, \beta \in \bN^n$ with  
$|\alpha|\ne |\beta|$ and $P(z)=a z^\alpha+ b z^\beta$
for some $a, b\in\bC$. Let $\Lambda=\Lambda(\p)$ 
with $\Lambda(\xi)$ homogeneous. Then Conjecture \ref{GVC} 
holds for $\Lambda$, $P(z)$ and any $g(z)\in \cz$.
\end{corol}
%
%


\vskip10mm


{\small  \sc A. van den Essen, Department of Mathematics, Radboud University Nijmegen, Postbus 9010,   
 6500 GL Nijmegen, The Netherlands.}

{\em E-mail}: essen@math.ru.nl \\

{\small \sc R. Willems, Department of Mathematics, Radboud University Nijmegen, Postbus 9010, 6500 GL Nijmegen, The Netherlands.}

{\em E-mail}: R.Willems@math.ru.nl\\

{\small \sc W. Zhao, Department of Mathematics, Illinois State University,
Normal, IL 61790-4520.}

{\em E-mail}: wzhao@ilstu.edu.

\end{document}